\documentclass{amsart}
\usepackage{amsmath, amssymb, euscript,  enumerate}

\newtheorem{thm}{Theorem}[section]
\newtheorem{cor}[thm]{Corollary}
\newtheorem{lem}[thm]{Lemma}
\newtheorem{claim}[thm]{Claim}
\newtheorem{prop}[thm]{Proposition}

\theoremstyle{remark}
\newtheorem{rem}[thm]{Remark}
\numberwithin{equation}{section}

\newcommand{\R}{\mathbb R}

\newcommand{\tv}{\tilde v}

\newcommand{\dist}{\mathrm{dist}}
\newcommand{\supp}{\mathrm{supp}}

\def\XXint#1#2#3{{\setbox0=\hbox{$#1{#2#3}{\int}$}
     \vcenter{\hbox{$#2#3$}}\kern-.5\wd0}}

\begin{document}
\title{Type II extinction profile  of maximal Solutions to the Ricci flow
in $\R^2$}

\author{Panagiota Daskalopoulos$^*$}
\address{Department of Mathematics, Columbia University, New  York,
USA} \email{pdaskalo@math.columbia.edu}
\thanks{$*:$ Partially supported by the NSF grants DMS-01-02252, DMS-03-54639  and the EPSRC in the UK}
 \author{Natasa Sesum}
\address{Department of Mathematics, Columbia University, New York, USA}
\email{natasas@math.columbia.edu}

\begin{abstract} 
We consider the initial value   problem  
$u_t  = \Delta \log u$,   
$u(x,0) = u_0(x)\ge 0$ in $\R^2$, corresponding to the Ricci flow, namely  conformal evolution of the metric $u \, (dx_1^2 + dx_2^2)$ by 
Ricci curvature. It is well known that the 
maximal (complete) solution 
$u$ vanishes identically after time $T= \frac 1{4\pi} \int_{\R^2} u_0 $. Assuming 
that $u_0$ is compactly supported 
we describe precisely the Type II vanishing behavior   of $u$ at time $T$: we show the existence of an inner region with exponentially fast vanishing profile, which is, up to 
proper scaling,  a {\em soliton cigar solution}, 
and the existence of an outer region of persistence of a logarithmic cusp. This is the only Type II singularity which has been  shown to exist, so far,   in the Ricci Flow in any dimension. It recovers 
rigorously formal asymptotics  derived by J.R. King \cite{K}.
\end{abstract}
\maketitle

\section{Introduction}

We consider the Cauchy problem
\begin{equation}
\begin{cases} \label{eqn-u}
u_t  = \Delta \log u  & \mbox{in} \,\,  \R^2
\times (0,T)\\
u(x,0) = u_0(x) & x \in \R^2 ,
\end{cases}
\end{equation}
for the {\em logarithmic fast diffusion}  equation in $\R^2$,
with $T >0$ and  
initial data $u_0 $ non-negative,  bounded and compactly supported. 

It has been observed by S. Angenent and L. Wu \cite{W1, W2}
that equation \eqref{eqn-u} represents the evolution of the
conformally equivalent metric $g_{ij} = u\, dx_i\, dx_j$ 
under the {\it Ricci Flow}
 \begin{equation}\label{eqn-ricci}
\frac{\partial g_{ij}}{\partial t} = -2 \, R_{ij}
\end{equation}
 which evolves $g_{ij}$  by its Ricci curvature. 
The equivalence easily follows  from the observation that the
conformal metric $g_{ij} = u\, I_{ij}$ has scalar curvature $R = -
(\Delta \log u) /u$ and in two dimensions $R_{ij} = \frac 12 \,
R\, g_{ij}$. 

Equation \eqref{eqn-u} arises also in physical applications, as a
model for long Van-der-Wals interactions in thin  films of a fluid
spreading on a solid surface, if certain nonlinear fourth order
effects are neglected, see \cite{dG,B,BP}.

It is shown in
\cite{DD1} that given an   initial data $u_0 \geq 0$ with $\int_{\R^2} u_0 \, dx < \infty$ and a constant $\gamma  \geq 2 $, there exists a  solution
$u_\gamma$ of \eqref{eqn-u} with
\begin{equation}\label{eqn-intc}
\int_{\R^2} u_\lambda (x,t) \, dx = \int_{\R^2} u_0 \, dx  -
 2\pi \gamma \, t.
\end{equation}
The solution $u_\gamma$ exists up to the exact time $T=T_\gamma$,
which is determined in terms of the initial area and $\gamma$ by
$T_\gamma= \frac{1}{2\pi\, \gamma} \, \int_{\R} u_0 \, dx.$

We restrict our attention to   {\em maximal solutions} $u$ of
\eqref{eqn-u}, corresponding to the value $\gamma =2$ in
\eqref{eqn-intc}, which vanish at  time
\begin{equation}\label{eqn-mvt}
T= \frac{1}{4 \pi}  \, \int_{\R^2} u_0 (x) \, dx.
\end{equation}
It is shown  in \cite{DD1} and \cite{RVE}
that if $u_0$ is compactly supported, then the maximal solution $u$ which
extincts  at time $T$  satisfies the asymptotic behavior
\begin{equation}\label{eqn-agc}
u(x,t)  =  \frac {2  t}{|x|^2\log^2 |x|} \left ( 1 + o(1)
\right ), \qquad \mbox{as \,\, $|x| \to \infty$}, \quad 0 \leq t <
T.
\end{equation}
This bound, of course,  deteriorates as $t \to T$. Geometrically 
 \eqref{eqn-agc} 
corresponds to the condition that the conformal metric is
complete. The manifold can be visualized as a surface 
of finite area with an unbounded cusp.

J.R. King \cite{K} has formally analyzed the extinction
behavior of maximal solutions $u$ of \eqref{eqn-u}, as $t \to T^-$.
His analysis, for
radially symmetric and  compactly supported initial data, suggests the existence of  two
regions of different behavior:   in  the {\em outer region}
$(T-t)\, \log  r > T$ the "logarithmic cusp" exact solution  $2 t\,
/|x|^2 \, \log^2 |x|$ of equation  $u_t = \Delta \log u$ persists.
However, in the {\em inner region}   $(T-t)\, \log r \leq  T$ the
solution vanishes exponentially fast and approaches, after an
appropriate change of variables,   one of the soliton solutions
$U$ of equation $U_\tau = \Delta \log U$ on $-\infty < \tau <
\infty$  given by $U(x,\tau) = 1/( \lambda  |x|^2 + e^{4\lambda \tau})$,  with  $\tau=1/(T-t)$ and $\lambda=T/2$.

This behavior was established rigorously in the radially symmetric case
by the first author and M. del Pino in \cite{DD2}. The precise asymptotics
of the Ricci flow neckpinch in the compact case, on $S^n$, has been
established by Knopf and Angenent in \cite{AngKn}.

Our goal in this paper is to remove the assumption of radial
symmetry and establish  the vanishing  behavior of maximal solutions 
of \eqref{eqn-u} for  any non-negative compactly supported initial data.

To state the inner region behavior in a precise manner, we perform the change of variables
\begin{equation}\label{eqn-rbu}
\bar u(x,\tau) =  \tau^2 \, u(x, t), \qquad \tau = \frac 1{T-t}
\end{equation}
and
\begin{equation}\label{eqn-rtu}
\qquad \tilde u(y,\tau) = \alpha(\tau) \, \bar
u(\alpha(\tau)^{1/2} y,\tau),
\end{equation}
with
\begin{equation}\label{eqn-atau}
\alpha(\tau) = [\bar u(0,\tau)]^{-1} = [(T-t)^{-2} u(0,t)]^{-1}
\end{equation}
so that $\tilde u(0,\tau)=1$.

A direct computation shows that  the rescaled solution $\tilde u$
satisfies the equation

\begin{equation}\label{eqn-tu}
\tilde u_\tau = \Delta \log \tilde u + \frac {\alpha
'(\tau)}{2\alpha (\tau)} \, \nabla ( y \cdot \tilde u) +\frac{2
\tilde u}{\tau}.
\end{equation}
Then,  following result holds:
\smallskip

\begin{thm}\label{Mth1} (Inner behavior)
Assume that the initial data  $u_0$ is non-negative, bounded and compactly supported.  Then, 
\begin{equation}\label{eqn-lim}
\lim_{\tau \to \infty} \frac{ \alpha ' (\tau)}{2\, \alpha(\tau)} = T 
\end{equation}
and the 
rescaled solution $\tilde u$ defined by \eqref{eqn-rbu} -
\eqref{eqn-atau} converges,  uniformly on compact subsets of
$\R^2$, to the solution 
$$U(x) = \frac 1{  \frac T2 \, |y|^2 +
1}$$ of the steady state equation
$$\Delta \log  U + T\cdot   \nabla ( y \cdot  U)=0.$$
\end{thm}

\smallskip
Since for any maximal solution $T= (1/{4\pi}) \int_{\R^2} u_0(x)\, dx$, this theorem shows, in particular, that the limit of the rescaled solution
is uniquely determined by the area of the initial data. The uniqueness of
the limit has not previously shown in \cite{DD2} even under the assumption
of radial symmetry.  

To describe the vanishing  behavior of $u$ in the outer
region we first perform the cylindrical  change of variables
\begin{equation}\label{eqn-rv}
v(\zeta,\theta,t) = r^2\, u(r,\theta,t), \qquad \zeta =\log r
\end{equation}
with $(r,\theta)$ denoting the polar coordinates.  Equation  $u_t=\Delta \log u$ 
in cylindrical coordinates takes the form 
\begin{equation}\label{eqn-v}
v_t = \Delta_c \log v
\end{equation}
with $\Delta_c$ denoting the Laplacian in cylindrical coordinates defined as
$$\Delta_c \log v=(\log v)_{\zeta\zeta} + (\log v)_{\theta\theta}.$$

We then perform a further scaling setting
\begin{equation}\label{eqn-rtv}
\tilde v(\xi ,\theta, \tau) =  \tau^2\,  v (\tau \xi,\theta,t), \qquad \tau=
\frac 1{T-t}.
\end{equation}
A direct computation shows that $\tilde v$ satisfies the equation
\begin{equation}\label{eqn-tv}
\tau \, \tilde  v_{\tau} = \frac 1{\tau} (\log \tilde v)_{\xi\xi} + \tau \, (\log \tilde v)_{\theta\theta}
+ \xi \, \tilde v_{\xi} + 2\tilde v.
\end{equation}
The extinction behavior of $u$ (or equivalently of $v$) in the
outer region $\xi \geq T$,  is described in the following
result.

\smallskip

\begin{thm}\label{Mth2} (Outer behavior).  Assume that the initial data
$u_0$ is non-negative, bounded and compactly supported. Then,  the rescaled solution $\tilde v$ defined
by \eqref{eqn-rtv} converges,  as $\tau \to \infty$,  to the $\theta-$independent
steady state solution $V(\xi)$ of equation \eqref{eqn-tv}   given by
\begin{equation}\label{eqn-dV}
V(\xi) =
\begin{cases}
\frac{ 2 T}{\xi^2}, \qquad  &\xi > T \\
0,  \qquad  &\xi < T. 
\end{cases}
\end{equation}
 Moreover, the convergence is uniform on
the set  $(-\infty, \xi^-]\times [0,2\pi]$,  for any $-\infty <   \xi^- < T$, and on compact subsets  of $(T, +\infty) \times [0,2\pi]$. \end{thm}

\smallskip

Under the assumption of radial symmetry this  result follows
from the work of the first author and del Pino \cite{DD2}.

\medskip
The proof of Theorems \ref{Mth1} and \ref{Mth2} rely on sharp estimates on the
geometric {\em width} $W$  and on the {\em maximum curvature}
$R_{\max}$  of  maximal solutions near their extinction time $T$
derived in \cite{DH} by the first author and R. Hamilton. In
particular, it is found in \cite{DH} that the maximum curvature is
proportional to $1/{(T-t)^2}$, which does not go along with the
natural scaling of the problem which would entail blow-up of order
$1/{(T-t)}$. One says that the vanishing behavior  is {\em of type II}. 
The proof also makes an extensive use of  the Harnack estimate on  the curvature $R = -\Delta \log u /u$ shown by Hamilton and Yau  in \cite{HY}. 
Although the result in \cite{HY} is shown only for a compact surface evolving by the Ricci flow, we shall observe in section \ref{sec-prelim} that the result remains valid in our case as well. 
Finally, let us remark  that  the proof of the inner-region behavior is based on the classification of eternal complete solutions of the 2-dimensional Ricci flow, recently shown 
by the authors in \cite{DS}.

\section{Preliminaries}\label{sec-prelim}

In this section we will collect some  preliminary results which
will be used throughout the rest of the paper. For the convenience
of the reader, we start  with a brief description of the geometric
estimates in \cite{DH} on which the proofs of Theorems \ref{Mth1}
and \ref{Mth2} rely upon.

\subsection{\large Geometric Estimates.}\label{sec-ge}
In \cite{DH} the first author and R. Hamilton established upper
and lower bounds on the geometric  width $W(t) $ of the
maximal solution $u$ of \eqref{eqn-u}
and on the maximum curvature  $R_{\max}(t)=
\max_{x \in \R^2} R(x,t)$, with $R= - (\Delta \log u) /u $.

Let $F:\R^2\to[0,\infty)$ denote a 
proper function $F$, such that    $F^{-1}(a)$ is compact for every
$a\in [0,\infty)$.  The width of $F$ is defined to be the supremum
of the lengths of the level curves of $F$, namely 
$w(F) = \sup_c L\{F=c\}.$
The width $w$ of the metric $g$, as introduced in \cite{DH},  is defined  to be the infimum
$$w(g) = \inf_F w(F).$$

The estimates  in \cite{DH} depend on the time to
collapse $T-t$. However, they do not scale in the usual way.
More precisely:  

\begin{thm}[Daskalopoulos and Hamilton \cite{DH}]
 \label{thm-DH1}
There exist positive constants $c$ and $C$ for which
\begin{equation}\label{eqn-w}
c \, (T-t) \leq  W(t) \leq C\, (T-t)
\end{equation}
and
\begin{equation}\label{eqn-c}
\frac{c}{(T-t)^2} \leq R_{\max}(t) \leq \frac{C}{(T-t)^2}
\end{equation}
for all $0< t < T$.
\end{thm}

\subsection{The Hamilton-Yau Harnack estimate.}

In \cite{HY} Hamilton and Yau established   a  Harnack estimate on the
curvature $R$ of a compact  surface  evolving by the Ricci flow, in the case where the curvature $R$ changes sign.
Since the proof in \cite{HY} uses only local quantities,  the result and its  proof 
can be carried over to the  complete, non-compact case.

\begin{thm}[Hamilton and  Yau \cite{HY}]
\label{thm-harnack}
For any constants $E$ and $L$ we can find positive constants $A, B, C, D$ such that for  any complete solution to the Ricci flow on $\R^2$ which   at the
initial time $t=0$ satisfies 
\begin{equation}
\label{equation-lower-R}
R \ge 1-E
\end{equation}
and
$$\frac{1}{R+E}\frac{\partial R}{\partial t} - \frac{|\nabla R|^2}{(R+E)^2} \ge -L$$
then,  for all $t\ge 0$ we have
\begin{equation}
\label{equation-harnack0}
\frac{1}{R+E}\frac{\partial R}{\partial t} - \frac{|\nabla R|^2}{(R+E)^2}
+ F(\frac{|\nabla R|^2}{(R+E)^2}, R+E) \ge 0
\end{equation}
where
$$F(X,Y) = A + \sqrt{2B(X+Y) + C} + D\log Y.$$
\end{thm}

Integrating the above estimate along paths we obtain:

\begin{prop}
Under the assumptions of Theorem \ref{Mth1}, there exist  uniform constants $E >0$ and $C_1, C_2 > 0$ so that for every $x_1, x_2 \in \R^2$ and $T/2 < t_1 < t_2$
we have 
\begin{equation}
\label{equation-harnack}
\frac{1}{\sqrt{R(x_1,t_1)+E}} \ge \frac{1}{\sqrt{R(x_2,t_2)+E}}- C_1(t_2-t_1) - C_2\frac{\dist_{t_1}^2(x_1,x_2)}{t_2-t_1}.
\end{equation}
\end{prop}

\begin{proof}
By the Aronson-Benil\'an inequality  $R \ge -1/t  \ge -2/T = 1-E$ for all $t\in [T/2,T)$.
Hence the estimate \eqref{equation-harnack0} and the lower curvature bound on $R$ give
\begin{eqnarray}
\label{equation-partial}
\frac{\partial R}{\partial t} &\ge& \frac{|\nabla R|^2}{R+E} - (2A+\sqrt{C})(R+E) - \sqrt{2B}|\nabla R|
\nonumber \\&-&\sqrt{2B}(R+E)\sqrt{R+E} - D\log(R+E)(R+E)\nonumber \\
&\ge& \frac{|\nabla R|^2}{R+E} - A_1(R+E)\sqrt{R+E} - \frac{1}{4}\frac{|\nabla R|^2}{R+E} \\
&=& \frac{3}{4}\frac{|\nabla R|^2}{R+E} - A_1(R+E)\sqrt{R+E}\nonumber.
\end{eqnarray}
Take any two points $x_1, x_2 \in \R^2$ and $T/2 \leq t_1 \le t_2 <T$ and let $\gamma$ be a curve
connecting $x_1$ and $x_2$, such that  $\gamma(t_1) = x_1$ and $\gamma(t_2) = x_2$. 
Since 
$$\frac{d}{dt}R(\gamma(t),t) = \frac{\partial R}{\partial t} + \langle \nabla R,\dot{\gamma}\rangle $$
using  also (\ref{equation-partial}) we find 
\begin{eqnarray*}
\frac{d}{dt}R(\gamma(t),t) &\ge& \frac{3}{4}\frac{|\nabla R|^2}{R+E} - A_1(R+E)\sqrt{R+E} - C|\dot{\gamma}|^2(R+E)
- \frac{1}{4}\frac{|\nabla R|^2}{R+E} \\
&\ge& -A_3(R+E)^{3/2}(1 + |\dot{\gamma}|^2).
\end{eqnarray*}
Integrating  the previous equation along the path $\gamma$, gives   
$$\frac{1}{\sqrt{R(x_1,t_1)+E}} \ge \frac{1}{\sqrt{R(x_2,t_2)+E}} - C(t_2-t_1) - \int_{t_1}^{t_2}|\dot{\gamma}|_{g(t)}^2dt.$$
Due to the bound $R \ge 1-E$ we have for $t\ge s$
$$|\dot{\gamma}|^2_{g(t)} \le (1-E) \, e^{s t}|\dot{\gamma}|^2_{g(s)}$$
and if we choose the curve $\gamma$  to be the minimal geodesic with respect to metric $g(t_1)$
connecting $x_1$ and $x_2$ we obtain
$$\frac{1}{\sqrt{R(x_1,t_1)+E}} \ge \frac{1}{\sqrt{R(x_2,t_2)+E}} - C_1(t_2-t_1) - C_2\frac{\dist_{t_1}^2(x_1,x_2)}{t_2-t_1}$$
as desired. 
\end{proof}

\subsection{Monotonicity of Solutions.}

Our solution $u(x,t)$ to (\ref{eqn-u}) has compactly supported initial data. The classical
argument based on reflection, due to Alexandrov and Serrin,  proves that such solutions enjoy the following monotonicity in the radial
direction:
\begin{lem}
\label{lemma-monotonicity}
Under the assumptions of Therorem \ref{Mth1}, if  $\supp \, u_0(\cdot) \subset B_\rho(0)$, then 
\begin{equation}
\label{equation-monotonicity}
u(x,t) \ge u(y,t),
\end{equation}
for all $t\in (0,T)$ and every pair of points $x,y\in\R^2$ such that $|y| \ge |x|+\rho$.
\end{lem}

The proof of Lemma \ref{lemma-monotonicity} is the same as the proof of Proposition $2.1$ in \cite{AC}. For the reader's convenience we will briefly sketch it. 

\begin{proof}[Sketch of the proof]
Assume that $\rho=1$.  By the comparison principle for maximal solutions it easily follows that
if $K = \supp u(\cdot,0)$ and $K \subset \{x\in \R^2:\,\,x_2 > 0\}$,
then $u(x_1,x_2,t) \ge u(x_1,-x_2,t)$ for $x_1\in \R^+$ and $t\in
[0,T)$. Fix $x^0\in B_1$ and $x^1\in \partial B_{1+\delta}$ for
$\delta > 0$. Let $\Pi$ be a hyperplane of points in $\R^2$ which are equidistant
from $x^0$ and $x^1$. Then,  it easily follows
$$\dist(\Pi,\{0\}) \ge 1$$
which implies $x^0$ and $\supp u(\cdot,0)$ are in the same half-space with respect to 
$\Pi$. Since $x^1$ is the reflection of $x^0$ in $\Pi$, it follows
$u(x^0,t) \ge u(x^1,t)$. We can now let $\delta\to 0$ to get the claim.
\end{proof}

Notice that due to (\ref{eqn-agc}), for every $t\in (0,T)$ we can define $x_t$
to be such that $u(x_t,t) = \max_{\R^2}u(\cdot,t)$. An easy consequence of
Lemma \ref{lemma-monotonicity} is the following result about $\{x_t\}_{t\in(0,T)}$.

\begin{cor}
\label{cor-maximums}
For every $t\in (0,T)$, $x_t\in B_{2\rho}(0)$.
\end{cor}  

\section{Inner Region Convergence}\label{sec-irc}

This section is devoted to the proof of the inner region
convergence, Theorem \ref{Mth1} stated in the Introduction. We
assume, throughout this section, that $u$ is a smooth, maximal solution of \eqref{eqn-u} 
with compactly supported initial data $u_0$ and $u$ a maximal solution that vanishes at time $$T= \frac 1{4\pi} \int_{\R^2} u_0 \, dx.$$

\subsection{\large Scaling and convergence}\label{sec-sc}
We introduce a new scaling on the solution $u$ namely 
\begin{equation}\label{eqn-rbu2}
\bar u(x,\tau) = \tau^2 \, u(x,t), \qquad \tau = \frac 1{T-t},
\quad \tau \in (1/T, \infty).
\end{equation}
Then $\bar u$ satisfies the equation
\begin{equation}\label{eqn-bu}
\bar u_\tau = \Delta \log \bar u + \frac{2\bar u}{\tau}, \qquad
\mbox{on \,\, $1/T \leq \tau < \infty.$}
\end{equation}
Notice that under this transformation,   $ \bar R := - \Delta \log
\bar u/ \bar u$  satisfies the estimate
\begin{equation}\label{eqn-c1}
 \bar R_{\max}(\tau)  \leq C
\end{equation}
for some  constant $C < \infty$. This is a direct
consequence of Theorem \ref{thm-DH1}, since $\bar R_{\max} (\tau)
= (T-t)^2\,  R_{\max}(t)$.

\smallskip

For an increasing  sequence $\tau_k \to \infty$ we set
\begin{equation}\label{eqn-ruk}
\bar u_k(y,\tau) = \alpha_k \, \bar u(\alpha_k^{1/2}\, y, \tau
+\tau_k), \qquad (y,\tau) \in \R^2 \times (- \tau_k + 1/T, \infty)
\end{equation}
where
$$\alpha_k = [\bar u(0,\tau_k)]^{-1}$$
so that $\bar u_k(0,0)=1$,  for all $k$. Then, $\bar u_k$
satisfies the equation
\begin{equation}\label{eqn-uk}
\bar u_\tau = \Delta \log \bar u + \frac{2\bar u}{\tau+\tau_k}.
\end{equation}
Let
$$\bar R_k := - \frac{\Delta \log \bar u_k}{\bar u_k}.$$
Then, by \eqref{eqn-c1}, we have
\begin{equation}\label{eqn-Rkk}
 \max_{y \in \R^2} \bar R_k(y,\tau) \leq   C, \qquad
 -\tau_k + 1/T < \tau < +\infty. 
\end{equation}
We will  also derive a global  bound from bellow on $\bar R_k$.
The Aronson-Benil\'an inequality $u_t \leq u/t$, on  $ 0 \leq t <
T$, 
 gives the bound
$ R(x,t) \geq - 1/t$  on $ 0 \leq t < T$. In particular, $ R(x,t)
\geq - C$  on $ T/2 \leq t < T$, which in the new time variable
$\tau=1/(T-t)$ implies the bound
$$\bar R(x,\tau) \geq - \frac{C}{ \tau^2}, \qquad 2/T < \tau <
\infty.$$ Hence
$$\bar R_k(y,\tau) \geq - \frac C{(\tau+\tau_k)^2}, \qquad  -\tau_k +
2/T < \tau < +\infty.$$ Combining the above inequalities we get
\begin{equation}\label{eqn-Rk}
- \frac{C}{(\tau+\tau_k)^2} \leq \bar R_k(y,\tau) \leq C, \qquad
\forall (y,\tau) \in \R^2 \times (-\tau_k + 2/T, +\infty).
\end{equation}

Based on the above  estimates we will  now  show  the following
convergence result.

\begin{lem}\label{lem-ick}
For each sequence $\tau_k \to \infty$, there  exists a subsequence
$\tau_{k_l}$ of $\tau_k$,  for which the rescaled solution $\bar
u_{\tau_{k_l}}$ defined by \eqref{eqn-ruk} converges, uniformly on
compact subsets  of $\R^2 \times \R$,    to an eternal solution
$U$ of equation $U_\tau = \Delta \log U$  on $\R^2 \times \R$ with
uniformly bounded curvature and uniformly bounded width.  Moreover, the convergence is in
$C^\infty(K)$, for any $K \subset \R^2 \times \R$  compact.
\end{lem}

\begin{proof} 
Denote by $x_k = x_{t_k}$ the maximum point of $u(\cdot,t_k)$. 
First, instead of rescaling our
solution by $\alpha_k$ we can rescale it by $\beta_k =
[\bar{u}(x_k,\tau_k)]^{-1}$, that is, consider
$$\tilde{u}_k(y,\tau) =
\beta_k\bar{u}(\beta_k^{1/2}\, y,\tau+\tau_k), \qquad 
\tau \in (-\tau_k+1/T,\infty).$$  For $y_k = \beta_k^{-1/2}\, x_k$
we have $\tilde{u}_k(y_k,0) = 1$ 
and $\tilde{u}_k(\cdot,0) \le 1$ since $x_k$ is the maximum point
of $u(\cdot,t_k)$. Notice that $|y_k| \leq  2\rho \beta_k^{-1/2}$,  because 
$x_k \in B_{2\rho}$, by Corollary \ref{cor-maximums}. Since $\tilde u_k$ satisfies \eqref{eqn-uk}, standard arguments imply  that $\tilde{u}_k$ is uniformly bounded from above and below
away from zero on any compact subset of $\R^2\times \R$. In particular, there
are uniform constants $C_1 >0$ and $C_2 < \infty$ so that
\begin{equation}
\label{equation-quotient}
C_1 \le \frac{\alpha_k}{\beta_k} \le C_2.
\end{equation}
Let $K\subset \R^2$ be a compact set. By  
(\ref{equation-quotient}), for every compact set $K$ there is a compact set $K'$ so
that for all $y \in K$ we have $y \, (\frac{\alpha_k}{\beta_k})^{1/2} \in K'$,
for all $k$. Also, by  the previous estimates we have
$$C_1(K') \le \frac{\bar{u}(\beta_k^{1/2}\,z,\tau_k+\tau)}{\bar{u}(x_k,\tau_k)} = \tilde u_k(z,\tau) \le C_2(K')$$
for all $z\in K'$ and $\tau$ belonging to a compact subset of $
(-\infty ,\infty)$. Therefore, using \eqref{equation-quotient} and
remembering that $\alpha_k = [\bar u(0,\tau_k)]^{-1}$ we find
$$\bar u_k(y,\tau)=\frac{\bar{u}(\alpha_k^{1/2}y,\tau+ \tau_k)}{\bar{u}(0,\tau_k)} \le
\frac{1}{C_1}\frac{\bar{u}(\beta_k^{1/2}[(\frac{\alpha_k}{\beta_k})^{1/2}y],\tau_k+\tau)}{\bar{u}(x_k,\tau_k)}
\le \frac{C_2(K')}{C_1} = C_2(K).$$  
Similarly,
$$C_1(K) = \frac{C_1(K')}{C_2} \le \frac{\bar{u}(\alpha_k^{1/2} \, y,\tau_k+\tau)}{\bar{u}(0,\tau_k)}= \bar u_k (y,\tau_k)$$
for $y \in K$ and $\tau$ belonging to a compact set.
Hence, by the  classical
regularity theory  the sequence $\{ \bar u_k \}$ is equicontinuous
on compact subsets of  $\R^2 \times \R$. It follows that there
exists a subsequence $\tau_{k_l}$ of $\tau_k$ such that $\bar
u_{k_l} \to U$ on compact subsets of $ \R^2 \times \R$, where $U$
is an eternal  solution of equation
\begin{equation}\label{eqn-U}
U_\tau = \Delta \log U, \qquad \mbox{on}\,\, \R^2 \times \R
\end{equation}
with infinite area  $\int_{\R^2} U(y,\tau)\, dy  = \infty$ (since
$\int_{\R^2} \bar u_k(y,\tau)  = 2(\tau + \tau_k)$). In
addition the classical regularity theory of quasilinear parabolic
equations implies that $\{u_{k_l} \}$ can be chosen so that
$u_{k_l} \to U$ in $C^\infty(K)$, for any compact set $K  \subset
\R^2 \times \R$, with $U(0,0) = 1$.

It then follows that $\bar R_{k_l} \to \bar  R:= -( \Delta \log
U)/U$. Taking the limit $k_l \to \infty$  on both sides of
\eqref{eqn-Rk} we obtain the bounds
\begin{equation}\label{eqn-cU}
0 \leq \bar  R \leq C, \qquad \mbox{on \,\, $\R^2 \times
\R$.}
\end{equation}
Finally, to show that $U$ has uniformly bounded width,  we take the  limit $k_l \to \infty$ in \eqref{eqn-w}. 
\end{proof}

\smallskip

As direct consequence of Lemma  \ref{lem-ick} and the classification result
of eternal solutions to the complete Ricci flow on $\R^2$, recently showed in \cite{DS}, we obtain  the following convergence result.

\begin{thm}\label{thm-ick}
For each sequence $\tau_k \to \infty$, there  exists a subsequence
$\tau_{k_l}$ of $\tau_k$ and numbers $\lambda, \bar{\lambda} >0$  for which the
rescaled solution $\bar u_{\tau_{k_l}}$ defined by \eqref{eqn-ruk}
converges, uniformly on compact subsets of $\R^2 \times \R$,  to the
soliton solution $U$ of the Ricci Flow given by
\begin{equation}\label{eqn-soliton}
U(y,\tau) = \frac 1{\lambda \, |y|^2 + e^{4 \bar{\lambda} \tau}}.
\end{equation}
Moreover, the convergence
is in $C^\infty(K)$,  for any $K \subset \R^2 \times \R$, compact.
\end{thm}

\begin{proof}
 
>From Lemma \ref{lem-ick}, $\bar u_{\tau_{k_l}} \to U$,
where $U$
 is an eternal solution of 
$U_t= \Delta\log U$, on $\R^2 \times \R$,  with
uniformly bounded width, such that $\sup_{\R^2}R(\cdot,t) \le C(t) <
\infty$ for every $t\in (-\infty,\infty)$. The main result in \cite{DS}
shows that the limiting solution $U$ is a
soliton of the form $U(x,\tau) = \frac 2{\beta \, (|x-x_0|^2 + \delta \,
e^{2\beta t})}$, with $\beta >0$, $\delta >0$,
 which under the
condition $U(0,0) =1$ takes the form $U(x,\tau) = \frac{1}{\lambda|x-x_0|^2 + e^{4\bar{\lambda}\tau}}$, with
$\lambda, \bar{\lambda} >0$.

\smallskip

It remains to show that the   limit $U(\cdot,\tau)$ is rotationally symmetric around the origin, that is, $x_0 = 0$.
This  will   follow  from Lemma \ref{lemma-monotonicity} and Lemma \ref{lem-ick}.  Notice that $\lim_{k\to\infty}\alpha_k = \infty$.
Since $\bar{u}_k(\cdot,\tau_k)$ converges uniformly on compact
subsets of $\R^2\times \R$ to a cigar soliton $U(y,0)$,
we have that
\begin{eqnarray*}
\bar{u}(0,\tau_k) &=& \tau_k^2u(0,t_k)\approx  \frac{1}{\lambda|x_0|^2 + e^{4\bar{\lambda}\tau_k}} \\
&\le& e^{-4\bar{\lambda}\tau_k} \to 0,
\end{eqnarray*}
as $k\to\infty$ and therefore $\lim_{k\to\infty}\alpha_k = \infty$. 
Lets us express  $\bar u=\bar u (r,\theta,\tau)$ in polar coordinates.
For every $r >0 $ there is $k_0$ so that $\alpha_k^{1/2}\, r > 1$ for
$k\ge k_0$. By Lemma \ref{lemma-monotonicity}
\begin{eqnarray*}
\min_{\theta}\bar{u}(\alpha_k^{1/2} \, r,\theta,\tau_k) &\ge& \max_{\theta} 
\bar{u}(\alpha_k^{1/2}\,  r+1,\theta,\tau_k) \\
&=& \max_{\theta}\bar{u}(\alpha_k^{1/2}(r+\alpha_k^{-1/2}),\theta,\tau_k)
\end{eqnarray*}
which implies
$$\min_{\theta}\bar{u}_k(r,\theta,0) \ge \max_{\theta}\bar{u}_k(r+\alpha_k^{-1/2},\theta,0).$$
Let $k\to\infty$ to obtain
$$\min_{\theta}U(r,\theta,0) \ge \max_{\theta} U(r,\theta,0)$$
which yields the limit $U(r,\theta,0)$ is radially symmetric with respect to
the origin and therefore $x_0=0$, implying that  $U$ is of  the form (\ref{eqn-soliton}). 
\end{proof}

\subsection{Further behavior} \label{sec-fb} We will now use the
geometric properties of the rescaled solutions and their limit, to
further analyze their vanishing behavior. Our analysis will be similar
to that in \cite{DD2}, applicable to the nonradial case as well.
However, the uniqueness of the limit along sequences $\tau_k \to \infty$ which will be shown in Theorem \ref{thm-curvature-limit}, 
is an improvement of the results in \cite{DD2}, even in the radial
case.

We begin by observing that rescaling back in the original $(x,t)$
variables, Theorem \ref{thm-ick} gives   the following asymptotic
behavior of the maximal solution  $u$ of \eqref{eqn-u}.

\begin{lem}\label{cor-ick1}
Assuming that along a sequence $t_k \to T$, the sequence $\bar
u_k$ defined by \eqref{eqn-ruk} with $\tau_k = (T-t_k)^{-1}$
converges to the soliton solution $U_\lambda$,  on compact subsets
of $\R^2 \times \R$, then along the sequence $t_k$ the  solution
$u(x,t)$ of \eqref{eqn-u} satisfies the asymptotics
\begin{equation}\label{eqn-asu1}
 u(x,t_k) \approx \frac{(T-t_k)^2} {\lambda \, |x|^2 + \alpha_k},
 \qquad \mbox{on} \quad |x| \leq \alpha_k^{1/2} \, M
\end{equation}
for  all $M >0$. In addition, the curvature $R(0, t_k) = - \Delta
\log u(0,t_k)/u(0,t_k)$ satisfies
\begin{equation}\label{eqn-limc}
\lim_{t_k \to T} (T-t_k)^2 \, R(0,t_k) = 4\, \lambda.
\end{equation}
\end{lem}

The proof of Lemma above is the same as the proof of Lemma $3.3$ in
\cite{DD2}. The following Lemma provides a sharp bound from below on the
 maximum curvature $4\, \lambda$ of the limiting solitons.

\begin{lem}\label{lem-bl} Under the assumptions of  Theorem \ref{Mth1}
the constant $\lambda$  in each   limiting solution
\eqref{eqn-soliton} satisfies
$$\lambda  \geq \frac{T}{2}.$$
\end{lem}

\begin{proof} We are going to use the estimate proven in Section 2 of
\cite{DH}. It is shown there   that if at time $t$ the  solution $u$
of \eqref{eqn-u} satisfies the scalar curvature bound $R(t) \geq -
2\, k(t)$, then the width $W(t)$ of the metric $u(t)\,
(dx_1^2+dx_2^2)$ (c.f. in Section \ref{sec-ge} for the definition)
satisfies the bound
$$
W(t)  \leq  \sqrt{k(t)} \, A(t) =   4  \pi \, \sqrt{k(t)} \, (T-t). 
$$
Here $A(t) = 4 \pi (T-t)$ denotes the area of the plane with
 respect
to the conformal metric $u(t)\, (dx_1^2+dx_2^2)$.
Introducing polar coordinates $(r,\theta)$, let  
$$\bar{U}(r,t)
= \max_{\theta}u(r,\theta,t) \quad \mbox{and} \quad \underbar{U}(r,t) =\min_{\theta}u(r,\theta,t).$$ Then
$$\underbar{U}(r,t) \le u(r,\theta,t) \le \bar{U}(r,t)$$
implying the bound 
\begin{equation}
\label{eqn-111}
W(\underbar{U}(t)) \le W(t) \le 4  \pi \, \sqrt{k(t)} \, (T-t).
\end{equation}
Observe next that the Aronson-Benil\'an inequality on $u$ implies the
bound
$R(x,t) \ge -{1}/{t}.$
Hence we can take $k(t) = \frac{1}{2t}$ in (\ref{eqn-111}). Observing
that for  the  radially symmetric solution $\underbar{U}$ the width
$W(\underbar{U}) = \max_{r\ge 0}2\pi r\sqrt{\underbar{U}}(r,t)$, we
conclude the pointwise estimate
\begin{equation}
\label{equation-222}
2 \pi  r \sqrt{\underbar{U}}(r,t) \le \frac{4\pi(T-t)}{\sqrt{2t}}, \qquad r\ge 0, \,\, 0<t<T.
\end{equation}
By Lemma (\ref{lemma-monotonicity}), 
\begin{equation}
\label{equation-333}
r\sqrt{u}(r+\rho,\theta,t) \le r\sqrt{\underbar{U}}(r,t), \qquad \mbox{for} \,\, r > 0.
\end{equation}
For a sequence $t_k\to T$,  let $\alpha_k = [\bar u(0,\tau_k)]^{-1}$, $\tau_k = 1/(T-t_k)$, as before. 
Using (\ref{eqn-asu1}), (\ref{equation-222}) and (\ref{equation-333}) we find
$$\frac{r(T-t_k)}{\sqrt{\lambda (r+\rho)^2+\alpha_k}} \le \frac{2(T-t_k)}{\sqrt{2t_k}}, \qquad
r\le M\alpha_k^{1/2},$$
for any positive number $M$. Hence, when $r = M\,\alpha_k^{1/2}$ 
we obtain the estimate
$$ \frac{M\, \alpha_k^{1/2} }{\sqrt{\lambda (M+\rho \, \alpha_k^{-1/2})^2 \, \alpha_k  +
\alpha_k}}   \leq \frac{2\, }{\sqrt{2\,t_k}}$$ or
$$\frac{M}{\sqrt{\lambda \, (M+\rho \, \alpha_k^{-1/2})^2 + 1}}   \leq
\frac{2 }{\sqrt{2t_k}}.$$ Letting $t_k \to T$ and taking
squares on both sides,  we obtain
$$\frac{1}{\lambda  + 1/M^2}   \leq
\frac{2}{T}.$$ Since $M >0$ is an arbitrary number, we finally
conclude $\lambda \geq T/2$, as desired.
\end{proof}

We will next provide a bound on the behavior of $\alpha(\tau)=
\tau^2 \, \bar u(0,\tau)$,  as $\tau \to \infty.$ In particular,
we will prove \eqref{eqn-lim}. 
We begin by a simple consequence of Lemma \ref{lem-bl}.

\begin{lem}\label{lem-altau1} Under the assumptions of Theorem
\ref{thm-ick} we have
\begin{equation}\label{eqn-altau}
 \liminf_{\tau \to \infty} \frac{\alpha'(\tau)}{\alpha(\tau)} \geq
4\, \lambda_0
\end{equation}
with $\lambda_0 = T/2$.
\end{lem}

The proof of Lemma \ref{lem-altau1} is the same as the proof of
Lemma $3.5$ in \cite{DD2}.

\begin{cor}\label{cor5}
Under the hypotheses of Theorem \ref{thm-ick},  we have
\begin{equation}\label{eqn-asa}
\alpha(\tau) \geq  e^{2T  \tau  + o(\tau)}, \qquad
\mbox{as \,\, $\tau \to \infty$}.
\end{equation}
\end{cor}

The next Proposition will be crucial in establishing the outer
region behavior of $u$.

\begin{prop}\label{prop-1}  Under the hypotheses of Theorem
\ref{Mth1}, we have
\begin{equation}\label{eqn-asa4}
\lim_{\tau \to \infty} \frac{\log \alpha(\tau)}{\tau}  = 2T.
\end{equation}
\end{prop}

\begin{proof} 
See Proposition $3.7$ in \cite{DD2}.
\end{proof}

A consequence of Lemma \ref{cor-ick1} and Proposition \ref{prop-1} is the following result, which will be used in the next section. 

\begin{cor}\label{cor-astv1}
 Under the assumptions of Lemma \ref{lem-ick} the rescaled solution
$\tilde v$ defined by \eqref{eqn-rtv} satisfies
$$\lim_{\tau \to \infty} \tilde v(\xi,\theta,\tau) =0, \qquad \mbox{uniformly
on} \,\, (\xi,\theta) \in (-\infty, \xi^-] \times [0,2\pi]$$ for all $\xi^- < T$.
\end{cor}
So far we have showed that $\bar{\lambda} = \lim_{\tau\to\infty}\frac{\log\alpha(\tau)}{\tau} = T/2$
and that $\lambda \ge T/2$. In the next theorem we will show that actually $\lambda = T/2$.
Theorem \ref{thm-curvature-limit} is an  improvement of the results  in \cite{DD2},  
since it leads to the uniqueness of a cigar soliton limit.

\begin{thm}
\label{thm-curvature-limit}
$\lim_{\tau\to\infty}\bar{R}(0,\tau) = 2T$.
\end{thm}

We will first prove the following lemma.

\begin{lem}
\label{lemma-xi}
For every $\beta > 1$ and for every sequence $\tau_i\to\infty$ there is a sequence
$s_i\in (\tau_i,\beta\tau_i)$ such that $\lim_{i\to\infty}\bar{R}(0,s_i) = 2T$.
\end{lem}

\begin{proof}
By definition 
$$(\log\alpha(\tau))_{\tau} = \bar{R}(0,\tau) - \frac{2}{\tau}.$$
Therefore
$$\log\alpha(\beta\tau_i) - \log\alpha(\tau_i) = (\bar{R}(0,s_i)-\frac{2}{s_i})(\beta - 1)\tau_i.$$
Since $\log \alpha(\tau) = 2T\tau + o(\tau)$, by Proposition \ref{prop-1},  we conclude 
$$(\bar{R}(0,s_i) - \frac 2{s_i})(\beta - 1) = 
\frac{(2T\beta\tau_i + o(\beta\tau_i)) - (2T\tau_i + o(\tau_i))}{\tau_i}$$
which yields
$$\bar{R}(0,s_i) = 2T + \frac 2{s_i} + \frac{o(\beta\tau_i) + o(\tau_i)}{\tau_i}$$
readily implying the Lemma. 
\end{proof}

\begin{proof}[Proof of Theorem \ref{thm-curvature-limit}]
By Lemma \ref{lem-bl},  we have  $\lambda \ge T/2$. 
Assume there is a sequence $\tau_i\to\infty$ such that $\lim_{i\to\infty}\bar{R}(0,\tau_i) = 4\lambda$,
where $4\lambda = 2T + \delta$ for some $\delta > 0$. We know that $\bar{R}(0,\tau) \le \tilde{C}$
for a uniform constant $\tilde{C}$. Choose $\beta > 1$ so that the following two conditions hold
\begin{equation}
\label{equation-first}
\frac{1}{\beta\sqrt{2\tilde{C}}} > C(\beta - 1)
\end{equation}
and
\begin{equation}
\label{equation-second}
2T + \frac \delta2 > 2T \left (\frac{\beta}{1-C(\beta - 1)\sqrt{2T}} \right )^2
\end{equation} 
for some uniform constant $C$ to be chosen   later. Notice that both  (\ref{equation-first})
and (\ref{equation-second}) are possible by choosing   $\beta > 1$suffciently  close to $1$.
By Lemma \ref{lemma-xi} find a sequence $s_i\in (\tau_i,\beta\tau_i)$ so that
$\lim_{i\to\infty}\bar{R}(0,s_i) = 2T$.
Let  $T/2 < t < T$. Then $R(x,t) \ge -\frac{2}{T} = 1-E$.
Hamilton-Yau Harnack estimate (\ref{equation-harnack}), applied to
$t_i$ (where $\tau_i = \frac{1}{T-t_i}$) 
and $\bar t_i>t_i$ (where $\frac{1}{T-\bar t_i} = \beta\tau_i$ for $\beta > 1$), 
yields
$$\frac{1}{\sqrt{\bar{R}(0,\tau_i)+ \frac{E}{\tau_i^2}}} \ge \
\frac{\tau_i}{s_i\sqrt{\bar{R}(0,s_i) + \frac{E}{s_i^2}}} - C\frac{s_i-\tau_i}{s_i}.$$ 
Notice that due to our choice of $\beta$ in (\ref{equation-first}) we have
\begin{eqnarray*}
\frac{\tau_i}{s_i\sqrt{\bar{R}(0,s_i) + \frac{E}{s_i^2}}} &\ge& \frac{\tau_i}{s_i\sqrt{2\tilde{C}}} \\
&\ge& \frac{1}{\beta\sqrt{2\tilde{C}}} \ge C(\beta - 1) \ge C\frac{s_i-\tau_i}{s_i}.
\end{eqnarray*}
Therefore 
\begin{eqnarray*}
\sqrt{\bar{R}(0,\tau_i) + \frac{E}{\tau_i^2}} &\le& \frac{1}{\frac{\tau_i}
{s_i\sqrt{\bar{R}(0,s_i) + \frac{E}{s_i^2}}} - C\frac{s_i-\tau_i}{s_i}} \\
&=& \frac{s_i\sqrt{\bar{R}(0,s_i)+\frac{E}{s_i^2}}}{\tau_i - C(s_i-\tau_i)
\sqrt{\bar{R}(0,s_i)+\frac{E}{s_i^2}}}.
\end{eqnarray*}
Denote by $A = \sqrt{\bar{R}(0,s_i) + \frac{E}{s_i^2}}$. Since the function
$f(x) = \frac{Ax}{\tau_i - CA(x-\tau_i)}$, for $x\in [\tau_i,\beta\tau_i]$
is increasing, we conclude 
\begin{eqnarray}
\label{equation-opposite}
\sqrt{\bar{R}(0,\tau_i) + \frac{E}{\tau_i^2}} &\le& \frac{\beta\tau_i\sqrt{\bar{R}(0,s_i)+\frac{E}{s_i^2}}}{\tau_i - 
C(\beta\tau_i-\tau_i)\sqrt{\bar{R}(0,s_i)+\frac{E}{s_i^2}}} \\
&=& \frac{\beta\sqrt{\bar{R}(0,s_i)+\frac{E}{s_i^2}}}{1 - 
C(\beta - 1)\sqrt{\bar{R}(0,s_i)+\frac{E}{s_i^2}}}. 
\end{eqnarray}
Letting  $i\to\infty$ in (\ref{equation-opposite}) we get 
$$\sqrt{4\lambda} \le \frac{\beta\sqrt{2T}}{1-C(\beta-1)\sqrt{2T}}$$
which implies 
$$2T + \delta = 2T \left (\frac{\beta}{{1-C(\beta-1)\sqrt{2T}}} \right )^2$$
contradicting  our choice of $\beta$ in (\ref{equation-second}).
\end{proof}

\subsection{Proof of Theorem \ref{Mth1}} We finish this section with
the proof of Theorem \ref{Mth1} which easily follows from the
results in Sections \ref{sec-sc} and \ref{sec-fb}. Take any 
sequence  $\tau_k \to \infty$. Observe that by Theorem
\ref{thm-curvature-limit} 
\begin{equation}
\label{eqn-altk}
\lim_{k \to \infty} \frac{\alpha'(\tau_k)}{\alpha(\tau_k)} = 2T.
\end{equation}
By the definitions of $\tilde u$ and $\bar u_k$ (\eqref{eqn-rtu} and
\eqref{eqn-ruk} respectively)   we have $\tilde u(y,\tau_k) =
\bar u_k(y,0)$.
By Theorem \ref{thm-ick}, we have $\bar u_k \to U_{\frac T2}$ and
therefore   $$\tilde
u(y,\tau_k) \to U_{\frac T2}(y,0) = \frac 1{\frac{T}{2} \, |y|^2
+1}.$$ The limit $U_{\frac T2}$ does not depend on the sequence $t_k\to T$ and  the proof of Theorem \ref{Mth1} is now complete.
\qed

\section{Outer Region Asymptotic Behavior}\label{sec-orc}

We assume,  throughout  this section,   that $u$ is a positive,
smooth, maximal solution of \eqref{eqn-u} satisfying the assumptions of Theorem
\ref{Mth2} which vanishes at time $$T=\frac 1{4\pi} \int_{\R^2} u_0 \, dx.$$
As in Introduction we consider the solution $ v(\zeta,\theta,t) = r^2\, u(r,\theta,t)$,
$\zeta =\log r$,  of the equation \eqref{eqn-v} in cylindrical coordinates. We
next set
\begin{equation}\label{eqn-rbv1}
\bar v(\zeta,\theta, \tau)  = \tau^2 \, v(\zeta,\theta, t), \qquad \tau = \frac 1{T-t}.
\end{equation}
 and
\begin{equation}\label{eqn-rtv1}
\tilde v(\xi ,\tau) =  \bar v (\tau \xi,\tau).
\end{equation}
The function $\tilde  v$ satisfies the equation 
\begin{equation}\label{eqn-tilv}
\tau \, \tilde  v_{\tau} = \frac 1{\tau} (\log \tilde v)_{\xi\xi} +\tau  (\log \tilde v)_{\theta\theta} + \xi \, \tilde v_{\xi} + 2\tilde v.
\end{equation}
Note  that the curvature $R=-\Delta_c \log v/v$,  is given in terms of $\tv$ by 
\begin{equation}
\label{eqn-tcurvature}
R(\tau \xi,\theta,t)  =- \frac{( \log \tv)_{\xi\xi} (\xi,\theta,t)+ \tau^2 ( \log \tv)_{\theta\theta}(\xi,\theta,t)}{\tv}. 
\end{equation}  
Moreover,  the
area of $\tilde v$ is  constant, in particular
\begin{equation}\label{eqn-mtv2}
\int_{-\infty}^\infty \int_0^{2\pi} \tilde v(\xi,\theta, \tau)\, d\theta \, d\xi =4\pi , \qquad
\quad\forall \tau.
\end{equation}

We shall show that, $\tilde v(\cdot, \tau)$ converges,  as $\tau
\to \infty$,  to  a  $\theta$-independent steady state of equation \eqref{eqn-tilv},
namely  to a solution of the linear first order equation
\begin{equation}\label{eqn-V2}
 \xi \, V_{\xi} + 2 V =0.
 \end{equation}
The area condition \eqref{eqn-mtv2} shall imply that
\begin{equation}\label{eqn-mV}
\int_{-\infty}^{\infty} V(\xi)\, d\xi =2.
\end{equation}
Positive solutions of equation \eqref{eqn-V2} are of the form
\begin{equation}\label{eqn-rV2}
V(\xi) = \frac{\eta }{\xi^2}
\end{equation}
where $\eta >0$ is any constant. These solutions become singular
at $\xi=0$ and in particular are non-integrable at $\xi=0$, so
that they do not satisfy the area condition \eqref{eqn-mV}.
However,  it follows from Corollary \ref{cor-astv1} that $V$ must
vanish in the interior region $ \xi <  T$. 
We will show that while  $\tilde v(\xi,\theta,\tau) \to 0$, as $\tau
\to \infty$ on  $(-\infty, T)$, we have  $\tilde v(\xi,\theta,\tau) \geq c
>0$, for $\xi > T$  and that actually $\tilde v(\xi,\theta,\tau)
\to 2\, T /\xi^2$,  on $(T, \infty)$, as stated in Theorem
\ref{Mth2}.

The rest of the section is devoted to the proof of Theorem
\ref{Mth2}.
We begin by showing the following  properties of the rescaled
solution $\tilde v$.

\begin{lem}\label{lem-ptv}  The rescaled solution $\tilde v$ given by
 \eqref{eqn-rbv1} - \eqref{eqn-rtv1} has the  following properties: \\
\noindent i. $\tilde v(\cdot,\tau) \leq C$,
for a constant $C$ independent of $\tau$.\\
\noindent ii. For any  $\xi^- < T$,   $ \tilde v(\cdot,\tau) \to 0$, as $\tau \to \infty$,
uniformly on $(-\infty,\xi^-] \times [0,2\pi]$. \\
\noindent iii. Let $\xi (\tau) = (\log \alpha(\tau))/2\tau$, with
$\alpha(\tau) = [\tau^2 \, u(0,t)]^{-1}.$ Then, there exists $\tau_0 >0$ and a constant $\eta >0$,
independent of $\tau$, such that
\begin{equation}\label{eqn-bbtilv}
 \tilde v(\xi,\theta,\tau) \geq \frac \eta{\xi^2}, \qquad {\mbox on}\,\,
\xi \geq  \xi (\tau),\,  \, \tau \geq \tau_0.
\end{equation}
In addition
\begin{equation}\label{eqn-xitau}
\xi(\tau) = T + o(1), \qquad \mbox{as} \,\, \tau \to \infty.
\end{equation}
\noindent iv. $\tilde v(\xi,\theta,\tau)$ also satisfies the upper bound
$$ \tilde v(\xi,\theta,\tau) \leq \frac {C}{\xi^2}, \qquad {\mbox on}\,\,
\xi >0, \,  \, \tau \geq \tau_0$$  for some constants $C >0$ and $\tau_0 >0$.

\end{lem}

\begin{proof} (i) One can easily show using the maximum principle that
$v(\zeta,\theta,t) \leq C/s^2$, for $s_0$ sufficiently large, with $C$ independent
of $t$. This implies the bound $\tilde v(\xi,\theta,\tau) \leq C/\xi^2$, for $\xi \, \tau > s_0$.
On the other hand,  by Corollary \ref{cor-astv1}, we have  
 $\tilde v (\xi,\theta,\tau) \leq C$, on $\xi < \xi^- <T$, with $C$ independent of $\tau$.
Combining the above, the desired estimate follows.

\noindent (ii) This is shown in Corollary \ref{cor-astv1}.

\noindent (iii) We have shown in the previous section that the rescaled solution $\bar u(x,\tau) = \tau^2 \, u(x,t)$, $\tau=1/(T-t)$,  defined by \eqref{eqn-rbu2} satisfies
the asymptotics $\bar u(x,\tau) \approx 1/(\frac T2 \, |x|^2 +
\alpha(\tau))$, when $|x| \leq   \sqrt{\alpha (\tau)}$. Hence
$$\tilde v (\xi(\tau), \theta, \tau) = \bar v (\xi(\tau)\, \tau,\theta,\tau) \approx 
\frac{e^{2\xi(\tau) \tau}}{\lambda \, e^{2\xi(\tau) \tau} + \alpha(\tau)} \approx \frac 1{\frac T2 +1}$$
if $\xi(\tau)=\frac{\log \alpha(\tau) }{2\tau}.$

Observe next that \eqref{eqn-xitau} readily follows from
Proposition \ref{prop-1}. Hence, it remains to show $\tilde v \geq
\eta /\xi^2$, for $\xi \in  [\xi(\tau), \infty)$, $\tau_0 \leq \tau < \infty$.
To this end, we will compare $\tilde v$ with the subsolution
$V_\eta(\xi,\theta) = {\eta}/{\xi^2}$ of equation \eqref{eqn-tilv}.
According to our claim above,  there exists  a  constant $\eta
>0$, so that
$$V_\eta(\xi(\tau),\theta) = \frac{\eta}{\xi(\tau)^2} \leq \tilde
v(\xi(\tau),\theta, \tau).$$ Moreover, by  the growth
condition \eqref{eqn-agc},  we can make
$$ \tilde v (\xi, \theta, \tau_0) > \frac {\eta}{\xi^2}, \qquad \mbox{on \,\, }
\xi \geq \xi (\tau_0)$$ by choosing $\tau_0 >0$ and $\eta$ sufficiently small. 
By the comparison principle,  \eqref{eqn-bbtilv} follows. 

\noindent (iv) Since $u_0$ is compactly supported and bounded, it follows that 
$u_0(r) \leq 2\, A/ (r^2\, \log^2 r)$, on $r>1$ for some $A >0$. Since  $2(t+A)/(r^2 \, \log^2 r)$ is an exact solution of equation \eqref{eqn-u}, 
it follows  by the comparison principle on  that  $u(r,t) \leq 2(t+A)/(r^2 \, \log^2 r$), 
for $r >1$, which readily
implies the  desired bound on $\tilde v$, with $C = 2(A+T)$.

\end{proof}

\begin{lem}
\label{lem-first-spherical}
For any compact set $K \subset (T,\infty)$,  there is a constant $C(K)$ for which
\begin{equation}
\label{equation-first-spherical}
\max_{\xi \in K} \left | \int_0^{2\pi}(\log\tilde{v})_{\xi} (\xi,\theta,\tau) \, d\theta \, \right | \le C(K), \qquad \forall \tau \ge 2/T.
\end{equation}
\end{lem}

\begin{proof}
We  integrate \eqref{eqn-tcurvature}   in $\theta$ variable and use the bounds  $R \ge - {1}/{t} \geq - 2/T$, for $t \geq T/2$,  and $\tv(\tau\xi,t) \le C$ shown in   Lemma \ref{lem-ptv},  to  get 
$$\int_0^{2\pi}  (\log\tilde{v} )_{\xi\xi}(\xi,\theta,\tau) \, d\theta \le C$$
for all $\tau = 1/(T-t) \geq 2/T$. 
We can now proceed as in the proof of Lemma $4.2$ in \cite{DD2}
to show (\ref{equation-first-spherical}).
\end{proof}

To simplify the notation, we set
$$R_c ( \zeta, \theta, t) = R(r, \theta,t), \qquad r=\log \zeta.$$

\begin{lem}\label{lem-inf}
For any compact set $K \subset (T,\infty)$,  there is a constant $C(K)$ such that for any $\xi_0 \in K$ and $\gamma  >0$
 $$\min_{[\xi_0,\xi_0 + \gamma \, \frac{\log\tau}{\tau}]\times[0,2\pi]} R_c(\xi \tau,\theta,t)\le \frac{C(K) \, \tau}{\gamma \, \log\tau}.$$
\end{lem}

\begin{proof}
Assume that for some $K$ and $\gamma$, $\min_{[\xi_0,\xi_0 + \gamma\, \frac{\log\tau}{\tau}] \times [0,2\pi]} R_c(\xi\tau,\theta,t) \ge \frac{M\tau}{\gamma \log \tau}$, 
for $M$ large. Then, it follows from \eqref{eqn-tcurvature} and the bound $\tilde v \leq C$ shown in   Lemma \ref{lem-ptv},  that
for every $\xi \in [\xi_0,\xi_0 + \gamma\, \frac{\log\tau}{\tau}]$ we have
\begin{eqnarray*}
\int_0^{2\pi}(\log\tilde{v})_{\xi\xi}(\xi,\theta,\tau)\, d\theta  &=& - \int_0^{2\pi}R_c(\xi\tau,\theta,t)) \, \tilde v(\tau\xi,\theta,t) \, d\theta \\
&\le& - \frac{C}{\xi^2}  \, \min_{[\xi_0,\xi_0 + \gamma \, \frac{\log\tau}{\tau}] \times[0,2\pi]} R_c(\xi\tau,\theta,t)  \\
&\le& -C_1(K) \, \frac{M \tau}{ \gamma \log \tau}
\end{eqnarray*}
which combined with  Lemma \ref{lem-first-spherical} implies
\begin{eqnarray*}
-C(K) &\le& \int_0^{2\pi} (\log\tilde{v})_\xi(\xi,\theta,\tau) \, d\theta \\
&\le& \int_0^{2\pi} (\log\tilde{v})_\xi (\xi_0,\theta,\tau) \, d\theta
- C_1(K) \, \frac{\gamma \log\tau}{\tau} \frac{M\tau}{\gamma \log \tau} \\
&=& \int_0^{2\pi} (\log\tilde{v})_\xi (\xi_0,\theta,\tau) \, d\theta - C_1(K)\, M \le C(K) - C_1(K) \, M 
\end{eqnarray*} 
impossible if $M$ is chosen sufficiently large.
\end{proof}  

\begin{prop}\label{prop-bcurv}
For every $K \subset (T,\infty)$ compact, there is a constant $C(K)$ depending only on $K$,  such that for any $\xi_0 \in K$ 
$$\max_{[\xi_0,\xi_0 + \frac{\log\tau}{\tau}] \times[0,2\pi]} R_c(\xi\tau,\theta,t)\le \frac{C(K)\, \tau}{\log\tau}.$$
\end{prop}

\begin{proof}
Let $\xi_1\in K$, $\theta\in [0,2\pi]$ and $\tau_1$ be arbitrary.
Choose $\xi_2$ such that $T < \xi_2 < \min K$  and $\tau_2 $  such  that $\xi_1\tau_1 = \xi_2\tau_2$. Since $\xi_2 < \xi_1$, then $\tau_2 > \tau_1$. Set $t_i= T -  1/{\tau_i}$, $i=1,2$. 
We next define  the set  
$A_{\xi_2} = \{\xi: \,\, \xi_2 \le  \xi \le \xi_2 + \gamma\frac{\log\tau_2}{\tau_2}\}$.  Let $\xi_0 \in
A_{\xi_2}$ and $\theta_2\in [0,2\pi]$ be such that
$$R_c(\xi_0\tau_2,\theta_2,t_2) = \min_{ (\xi,\theta) \in
A_{\xi_2}\times [0,2\pi] } \, R_c(\xi\tau_2,\theta,t_2)$$  and set $x_1 =
(e^{\xi_1\tau_1},\theta_1)$ and $x_2 = (e^{\xi_0\tau_2},\theta_2)$.
Since $\xi_1 \tau_1 = \xi_2 \tau_2 \leq \xi_0 \tau_2$,  then  $|x_1| \leq  |x_2|$.  
Denoting by $\dist_{t_1}(x_1,x_2)$ the distance with respect to
the metric $g_{t_1} = u(\cdot, t_1) \, (dx^2 + dy^2)$, we have:

\begin{claim}
For any $0 < \gamma <1$, there is a constant $C=C(K,\gamma)$ so that
$$\dist_{t_1}(x_1,x_2) \le \frac{C(K,\gamma)}{\tau_1^{1-\gamma}}.$$
\end{claim}

\noindent{\em Proof of Claim.}
We have seen in the proof of Lemma \ref{lem-ptv} that $u(x,t) \le \frac{C}{|x|^2\log^2|x|}$, 
for all $|x| \ge 1$ and all $t\in [0,T)$. If $\sigma$ is a euclidean geodesic with respect to $g_{t_1}=  u(\cdot, t_1) \, (dx^2 + dy^2)$,  connecting  $x_1$ and $x_2$, this implies
\begin{eqnarray}
\label{equation-dist}
\dist_{t_1}(x_1,x_2) &\le& \int_{\sigma} \sqrt{u}(\cdot,t_1)\,d\sigma \nonumber \\
&\le& C \, \frac{|e^{\xi_1\tau_1} - e^{\xi_0\tau_2}|}{e^{\xi_1\tau_1}\xi_1\tau_1} \nonumber \\
&\le& C\, \frac{e^{\xi_2\tau_2 + \gamma\log\tau_2} - e^{\xi_1\tau_1}}{e^{\xi_1\tau_1}\xi_1\tau_1} \nonumber \\
&\le& \frac{C}{\xi_1\tau_1}(e^{\gamma\log\tau_2} - 1) = \frac{C \, \xi_1^{\gamma -1}}{\xi_2^{\gamma}\tau_1^{1-\gamma}} - \frac C{\xi_1\tau_1} \nonumber \\
&\le& \frac{C \, \xi_1^{\gamma -1}}{\xi_2^{\gamma}\tau_1^{1-\gamma}} \le \frac{A(K,\gamma)}{\tau_1^{1-\gamma}}.
\end{eqnarray}  
\qed

To finish  the proof of the Proposition, we first apply the Harnack estimate  \eqref{equation-harnack} to obtain the inequality
$$\frac{1}{\sqrt{R_c(\xi_1\tau_1,\theta_1,t_1)+E}} \ge \frac{1}{\sqrt{R_c(\xi_0\tau_2,\theta_2,t_2)+E}} - C\, (t_2-t_1) - 
C\, \frac{\dist^2_{t_1}(x_1,x_2)}{t_2-t_1}.$$
By Lemma \ref{lem-inf}, using also that   $\xi_1\tau_1 = \xi_2\tau_2$ (since $\tau_2 > \tau_1$, we have $\xi_1 > \xi_2$),  we get 

\begin{eqnarray*}
\frac{1}{\sqrt{R_c(\xi_1\tau_1,\theta_1,t_1)+E}} &\ge& \frac{C(K,\gamma)\, \sqrt{\log\tau_2}}{\sqrt{\tau_2}} \, -  \frac{C\, (\xi_1 - \xi_2)}{\xi_2\tau_2} \, 
- \frac{C(K,\gamma)}{\tau_1^{1-2\gamma}} \, \frac {\tau_2} {(\tau_2-\tau_1)} \\
&=&\frac{C(K,\gamma)\, \sqrt{\log\tau_2}}{\sqrt{\tau_2}} \, - \frac{C_1(K,\gamma)}{\tau_2}
- \frac{C_2(K,\gamma)}{\tau_1^{1-2\gamma}}. \end{eqnarray*}
In the last inequality we used that 
$$\frac {\tau_2} {(\tau_2-\tau_1)} = \frac {1} {(1-\tau_1/\tau_2)}
= \frac {1} {(1-\xi_2/\xi_1)} = \frac {\xi_1} {(\xi_1-\xi_2)}$$
and that $(\xi_1-\xi_2)/\xi_2$ and $\xi_1/(\xi_1-\xi_2)$ depend only on the set $K$. 

Take $\gamma = \frac{1}{4}$. Using  that $\tau_2/\tau_1=\xi_1/\xi_2$
depends only on $K$, we  conclude  the inequalities 
\begin{eqnarray}
\label{equation-curv1}
\frac{1}{\sqrt{R_c(\xi_1\tau_1,\theta_1,t_1)+E}} &\ge& \frac{\tilde C (K) \, \sqrt{\log\tau_2}}{\sqrt{\tau_2}} \nonumber \\
&\ge& \frac{\tilde C_1 (K) \sqrt{\log\tau_1}}{\sqrt{\tau_1}}
\end{eqnarray}
for $\tau_1$ sufficiently large, depending only on  $K$. 
Estimate (\ref{equation-curv1}) yields the bound
$$R_c(\xi_1\tau_1,\theta_1,t_1) \le  \frac{C(K) \,\tau_1}{\log\tau_1}$$
 finishing the proof of the Proposition.  
\end{proof}

\begin{cor}\label{cor-error}  Under the assumptions of Theorem \ref{Mth2}, we have $$\lim_{\tau \to \infty} \frac 1{\tau} \int_0^{2\pi} (\log \tv)_{\xi\xi}(\xi,\theta,\tau) \, d\theta
=0$$
uniformly on compact subsets of $(T,\infty)$. 
\end{cor}
\begin{proof}
We begin by integrating  \eqref{eqn-tcurvature} in $\theta$ which gives 
$$ \int_0^{2\pi}  (\log \tv)_{\xi\xi}(\xi,\theta,\tau)\, d\theta = - \int_0^{2\pi} 
R_c(\xi  \tau,\theta,t)\, \tv(\xi,\theta,\tau)\,  d\theta, \quad  \tau=\frac1{T-t}.$$
Let $K \subset (T,\infty)$ compact. 
By the Aronson-Benil\'an inequality and Proposition \ref{prop-bcurv} $$ - \frac 1t \leq R_c(\xi\, \tau,\theta,t) \leq \frac{C(K) \, \tau}{\log \tau}.$$
Since $\tilde v \leq C$ (by   Lemma \ref{lem-ptv}),  we conclude 
\begin{equation}\label{eqn-000}
 \left | \frac 1\tau \int_0^{2\pi} (\log \tv)_{\xi\xi}(\xi,\theta,\tau)\, d\theta 
\right | \leq  \frac{C(K)}{\log \tau}
\end{equation}
from which the lemma directly follows. 
\end{proof}

We next  introduce the  new time variable
$$s = \log \tau = - \log (T-t), \qquad s \geq -\log T.$$ To simplify
the notation we still call $\tilde v(\xi,\theta, s)$ the solution $\tilde
v$  in the new time scale. Then, it is easy to compute that
$\tilde v(\xi,\theta, s)$ satisfies the equation
\begin{equation}\label{eqn-tilvs}
\tilde v_s = e^{-s} \, (\log \tilde v)_{\xi\xi} + e^s (\log \tilde v)_{\theta\theta}+ \xi\, \tilde v_\xi + 2\, \tilde v.
\end{equation}
For an increasing sequence of times $s_k \to \infty$, we let
$$\tilde  v_k (\xi,s) = \tilde v (\xi, s+s_k), \qquad - \log T-s_k  <
s < \infty.
$$
Each $\tilde  v_k$ satisfies the equation
\begin{equation}\label{eqn-vkk}
(\tilde  v_k)_s = e^{-(s+s_k)} (\log \tilde  v_k)_{\xi\xi} +  e^{s+s_k}(\log\tilde{v}_k)_{\theta\theta}
+ \xi \, (\tilde  v_k)_{\xi} + {2\tilde  v_k}
\end{equation}
and  the area condition
\begin{equation}\label{eqn-mcvk}
\int_{-\infty}^\infty\int_0^{2\pi}  \tilde v_k(\xi,\theta,s)\, d\theta \,  d\xi = 2 .
\end{equation}
Defining  the functions  
$$W_k(\eta,s) = \int_{\eta}^{\infty}
\int_0^{2\pi}\tilde{v}_k(\xi,\eta,s)\,d\theta\,d\xi, \quad \eta \in (T,\infty), \, 
 - \log T -s_k < s < \infty$$
we have:

\begin{lem}\label{prop-ock2}
Passing to a subsequence, $\{W_k\}$ converges uniformly on compact subsets of $\eta\in (T,\infty)$ to the time-independent steady state
${2T}/{\eta}$. In addition, for any $p \ge 1$ and $\xi_0\in (T,\infty)$, the solution $\tilde v_k(\xi,\theta,s)$ of \eqref{eqn-vkk} converges in $L^p([\xi_0,\infty)\times[0,2\pi])$ norm to ${2T}/{\xi^2}$.
\end{lem}

\begin{proof}
We first integrate \eqref{eqn-vkk} in $\theta$ and $\xi\in [\eta,\infty)$, for
$\eta\in (T,\infty)$,  to find that each $W_k$ satisfies the equation 
$$(W_k)_s = -\int_{\eta}^{\infty}\int_0^{2\pi}\frac{(\log \tilde v_k)_{\xi\xi}}{\tau_k(s)}\,d\theta\,d \xi 
+ \int_{\eta}^{\infty}\int_0^{2\pi} \xi \, (\tilde{v}_k)_{\xi}\,d\theta\,d\xi + 2 \, W_k$$
with $\tau_k(s)= e^{-(s+s_k)}$. 
Integrating  by parts the second term yields 
\begin{eqnarray*}
\int_{\eta}^{\infty}\int_0^{2\pi}\xi \, (\tilde{v}_k)_{\xi}\,d\theta\,d\xi &=&
 -W_k(\eta,s) - \eta\int_0^{2\pi}\tilde{v}_k(\eta,\theta,s)\,d\theta + \int_0^{2\pi}\lim_{\xi\to\infty}\xi \tilde{v}_k(\xi,\theta,s)\, d\theta \\
&=& -W_k(\eta,s) - \eta\int_0^{2\pi}\tilde{v}_k(\eta,\theta,s)\,d\theta 
\end{eqnarray*}
since due to our estimates on $\tilde{v}$ in Lemma \ref{lem-ptv},   we have 
$\lim_{\xi\to\infty} \xi \tilde{v}_k(\xi,\theta,s) = 0$, uniformly in $k$ and $\theta$.  
We conclude that 
$$(W_k)_s = -  \int_{\eta}^{\infty}\int_0^{2\pi}\frac{(\log \tilde v_k)_{\xi\xi}}{\tau_k(s)}\,d\theta\,d \xi \,
+ W_k + \eta \, (W_k)_{\eta}.$$
Let $K \subset (T,\infty)$ compact. Then,  by \eqref{eqn-000}
$$ \left  | \int_{\eta}^{\infty}\int_0^{2\pi}\frac{(\log \tilde v_k)_{\xi\xi}}{\tau_k(s)}\,d\theta\,d \xi \, \right | \leq \frac{C(K)}{s+s_k}.$$
Also,  by Lemma  \ref{lem-ptv} and  Proposition \ref{prop-bcurv},
there exists a constant $C=C(K)$ for which the bounds 
 \begin{equation}
\label{equation-uniform-est}
|W_k(\eta,s)| \le C, \quad |(W_k)_s(\eta,s)| \leq C,  \quad |(W_k)_{\eta}(\eta,s)| \le C
\end{equation}
hold, for $s\ge -\log T$. Hence, passing to a
subsequence, $W_k(\eta,s)$ converges uniformly on compact subsets of $(T,\infty)\times \R$
to a solution $W$ of the equation
$$W_s = \eta \,  W_{\eta} + W = (\eta \,  W)_{\eta} \qquad \mbox{on} \,\,\, (T,\infty)\times \R$$
with
$$\lim_{\eta\to T} W(\eta,s) = 2, \qquad s\in \R$$
and
$$\lim_{\eta\to\infty}W(\eta,s) = 0, \qquad s\in \R.$$
As in \cite{DD2},  one  can show that $W$ is completely determined by its boundary values at $T$, and it is is the steady state 
$$W(\eta,s) = \frac{2\, T}{\eta},  \qquad \eta > T, \,\, s \in \R.$$
To show the $L^p$ convergence, we first notice that by the
comparison principle 
$$v(\zeta,t) \le \frac{2\, T}{(\zeta - \zeta_0)^2}, \qquad \zeta \ge \zeta_0,\,\, 0 < t < T$$
for $\zeta_0=\log \rho$, with $\rho$ denoting the radius of the  support of $u_0$. 
This yields the bound
$$\tilde{v}_k(\xi,\theta,s) \le \frac{2\, T}{(\xi - \zeta_0/\tau_k(s))^2}, \qquad \xi \ge T.$$
By the triangle inequality  we have
\begin{eqnarray*}
\int_{\eta}^{\infty} \int_0^{2\pi} |\frac{2\, T}{\xi^2} - \tilde{v}_k |\, d\theta\,d\xi  &\le&
\int_{\eta}^{\infty} \int_0^{2\pi} \left ( \frac{2\, T}{(\xi - \zeta_0/\tau_k(s))^2} - \tilde{v}_k \right )\,d\theta\,d\xi \\
&+& \int_{\eta}^{\infty}\int_0^{2\pi} \left (\frac{2\, T}{(\xi - \zeta_0/\tau_k(s))^2} - \frac{2T}{\xi^2} \right ) \, d\theta\,d\xi
\end{eqnarray*}
where the second  integral converges to zero, as $k \to \infty$,  by the first part of the Lemma. It is easy to see
that the third integral converges  as well. This gives us the desired
$L^1$ convergence, which immediately implies the $L^p$ convergence,  since
$|\tilde{v}_k(\xi,\theta,s) - {2T}/{\xi^2}|$ is uniformly  bounded on $[\xi_0,\infty)\times [0,2\pi]$, 
for $\xi_0  \ge T$ and $s\ge -\log T$.
\end{proof}

\begin{rem}
\label{rem-pointwise}
The $L^p$ convergence in the previous Lemma, implies that there is a subsequence $k_l$ so that $\tilde{v}_{k_l} (\xi,\theta,s)
\to {2T}/{\xi^2}$ pointwise, almost everywhere on $(T,\infty)\times[0,2\pi]$. 
\end{rem}

Set $\tau_k(s) = e^{s+s_k}$. Since  
$$\frac{(\log\tilde{v})_{\xi\xi}(\xi, \theta,\tau)}{\tau}  + \tau \, (\log\tilde{v})_{\theta\theta}(\xi,
\theta,\tau) 
= -\tau \, R_c(\xi\, \tau,
\theta,\tau) \, \tilde v(\xi,\theta,t),$$
we can rewrite (\ref{eqn-vkk}) as 
\begin{equation}
\label{equation-rewrite}
(\tilde{v}_k)_s = -\frac{R_c }{\tau_k(s)}\, \tilde v_k  + \xi \, (\tilde  v_k)_{\xi} + {2 \, \tilde  v_k}.
\end{equation}
We divide the equation by $\tilde{v}_k$ and integrate it in $\theta$.  Denoting by 
$Z_k(\xi,s) = \int_0^{2\pi} \log \tilde v_k(\xi,\theta,s)d\theta$  we get 
$$(Z_k)_s = -\int_0^{2\pi}\frac{R_c}{\tau_k(s)} \, d\theta + \xi (Z_k)_{\xi} + 4\pi.$$
Notice that by Proposition (\ref{prop-bcurv}), we have 
\begin{equation}
\label{equation-curv-lim}
|\frac{R_c}{\tau_k(s)}| \le \frac{1}{\log\tau_k(s)}
\end{equation}
and that by Lemma \ref{lem-first-spherical}
$$|(Z_k)_{\xi}(\xi,s)| = \left |\int_0^{2\pi}(\log\tilde{v}_k)_{\xi}(\xi,\theta,s) \, d\theta
\right  | \le C(K)$$
for $\xi\in K$, a compact subset of $(T,\infty)$ and $s\ge -s_k - \log T$.
This also implies the bound
$$|(Z_k)_s(\xi,s)| \le C(K).$$

\begin{lem}\label{prop-ock3}
Passing to a subsequence, $\tilde Z_k(\xi,s)$ converges uniformly
on compact subsets of $(T,\infty) \times \R$
to a   solution  $Z$ of the  equation
\begin{equation}\label{eqn-V10}
Z_s = \xi \,   Z_{\xi} + 4\, \pi \qquad  \mbox{on} \,\,  (T, \infty) \times \R. 
\end{equation}
\end{lem}

\begin{proof}
Let $ E \subset (T, \infty) \times \R$
 compact. Then
according to the previous estimates, the sequence
 $\tilde Z_k$ is
equicontinuous on $E$, hence passing to a subsequence it converges
to a function $Z$. In addition, the estimate \eqref{equation-curv-lim} readily
implies that $Z$ is a solution of the first order equation
\eqref{eqn-V10}.
\end{proof}

\begin{lem}
\label{claim-right-thing}
The function $Z$ is given by 
\begin{equation}
\label{eqn-Z}
Z(\xi,s) = 2\pi \, \log \frac{2T}{\xi^2}, \qquad (\xi,s)\in  (T,\infty)\times \R.
\end{equation} \end{lem}

\begin{proof}
Since $\int_0^{2\pi}\log\tilde{v}_k(\xi,\theta,s) \, d\theta \to  Z(\xi,s)$,
uniformly in $\xi$ on compact subsets of $(T,\infty)$, then for any $A > 0$ we have 
$$\int_{\eta}^{\eta +A}\int_0^{2\pi}\log\tilde{v}_k(\xi,\theta,s)\,d\theta\,d\xi \to 
\int_{\eta}^{\eta+A} Z(\xi,s)\,d\xi.$$
By Remark \ref{rem-pointwise} and the dominated  convergence theorem it follows that for every $A > 0$ we have
$$\int_{\eta}^{\eta + A} 2\pi\, \log\frac{2T}{\xi^2}\, d\xi = \int_{\eta}^{\eta + A} Z(\xi,s)\,d\xi$$ 
implying that $Z$ is given by \eqref{eqn-Z}.

\end{proof}

We are finally in position to conclude the proof of Theorem \ref{Mth2}.

\smallskip

\noindent{\em Proof of Theorem \ref{Mth2}.}
We begin by observing that by Lemma \ref{lem-ptv}
$$\tilde v_k(\xi,\theta,\tau) \to 0, \qquad \mbox{as} \,\, \tau \to \infty$$ 
uniformly on $(-\infty,\xi^-] \times [0,2\pi]$, for any $-\infty < \xi^- < T$.

To show the convergence on the outer region, observe that by Lemma  \ref{prop-ock3} and  Lemma \ref{claim-right-thing} 
\begin{equation}
\label{equation-uniform}
\int_0^{2\pi}\log\tilde{v}_k(\xi,\theta,s) \, d\theta \to   2\pi\log\frac{2T}{\xi^2}
\end{equation}
uniformly on compact subsets  of $(T, \infty) \times (-\infty,\infty).$
Set 
$$\underline{v}_k(\xi,s) = \min_{\theta\in [0,2\pi]}\tilde{v}_k(\xi,\theta,s) \quad \mbox{and} \quad  \overline{v}_k(\xi,s) = \max_{\theta\in [0,2\pi]}\tilde{v}_k(\xi,\theta,s).$$ 
Let us recall that $u_0 \subset B_\rho(0)$. By the monotonicity property of the solutions
shown in  Lemma \ref{lemma-monotonicity}, we have
\begin{eqnarray}
\label{equation-squeeze}
2\pi \, \log \underline{v}_k(\xi,s) &\le& \int_0^{2\pi}\log\tilde{v}_k(\xi,\theta,s)\, d\theta   
\le 2\pi\, \log\overline {v}_k(\xi,s) \nonumber \\
&\le& 2\pi\, \log\frac{e^{2\xi\tau_k(s)}}{(e^{\xi\tau_k(s)}-1)^2} + 
2\pi \log\underline{v}_k(\xi + \frac{\log ( 1- \rho \, e^{-\xi\tau_k(s)})}{\tau_k(s)},s) \nonumber \\
&\le& 2\pi\log\frac{e^{2\xi\tau_k(s)}}{(e^{\xi\tau_k(s)}-\rho)^2} + 
\int_0^{2\pi}\log\tilde{v}_k(\xi + \frac{\log (1 - \rho \, e^{-\xi\tau_k(s)})}{\tau_k(s)},\theta,s)\, d\theta.
\end{eqnarray}
Combining (\ref{equation-uniform}) and (\ref{equation-squeeze}) yields
that 
$$\label{equation-radial-uniform}
2\pi\log\overline{v}_k(\xi,s) \to 2\pi\log\frac{2T}{\xi^2} \quad \mbox{and} \quad 2\pi\log \underline{v}_k(\xi,s) \to 2\pi\log\frac{2T}{\xi^2}$$
uniformly on compact subsets of $(T,\infty)\times \R$.
Since
$$\underline{v}_k(\xi,s) \le \tilde{v}_k(\xi,\theta, s) \le \bar{v}_k(\xi,s)$$
the above readily implies that 
$\tilde v_k(\xi,\theta,s) \to {2T}/{\xi^2}$
uniformly on compact subsets of $(T,\infty)\times [0,2\pi] \times (-\infty,\infty)$. Since the limit is independent of the sequence $s_k \to \infty$,  we conclude that $$\tilde v(\xi,\theta,\tau) \to \frac {2T}{\xi^2},
\qquad \mbox{as} \,\, \tau \to \infty$$
uniformly on compact subsets of $(T,\infty)\times [0,2\pi]$, which finishes 
the proof of Theorem \ref{Mth2}. \qed


\smallskip



\begin{thebibliography}{99}


\bibitem{Ang} Angenent, S.  The zero set of a solution of a parabolic
equation, {\em J. Reine Angew. Math.} 390 (1988), 79--96.

\bibitem{AngKn} Angenent, S., Knopf, D., Precise asymptotics of the Ricci flow 
neckpinch; preprint available on: www.ma.utexas.edu/~danknopf.

\bibitem{AB}
Aronson, D.G., B\'enilan P., R\'egularit\'e des solutions de
l'\'equation de milieux poreux dans ${\bf R}^n$, {\em C.R. Acad.
Sci. Paris, 288}, 1979, pp 103-105.

\bibitem{AC} Aronson, D.G., Caffarelli,L.A., The initial trace of a
solution of the porous medium equation, {\em Transactions of the
Amer. Math. Soc. 280} (1983), 351--366.

\bibitem{B} Bertozzi, A.L., The mathematics of moving contact lines in
thin liquid films, {\em  Notices Amer. Math. Soc. 45} (1998), no.
6, pp 689--697.

\bibitem{BP} Bertozzi, A.L., Pugh M.,  The lubrication approximation
for thin viscous films: regularity and long-time behavior of weak
solutions, {\em  Comm. Pure Appl. Math. 49} (1996), no. 2, pp
85--123.


\bibitem{CC}  Cao, H.-D., Chow, B.,  Recent developments on the Ricci
flow, {\em Bull. Amer. Math. Soc. (N.S.) 36} (1999), no. 1, pp
59--74.

\bibitem{Ch} Chow, B., The Ricci flow on the $2$-sphere. J.
Differential Geom. 33 (1991), no. 2, 325--334.

\bibitem{Ch2} Chow, B., On the entropy estimate for the Ricci flow on
compact $2$-orbifolds. J. Differential Geom. 33 (1991), no. 2,
597--600.


\bibitem{dG} de Gennes, P.G., Wetting: statics and dynamics,
{\em Reviews of Modern Physics, 57 No 3}, 1985, pp  827-863.

\bibitem{DD1} Daskalopoulos,P., del Pino M.A.,
On a Singular Diffusion Equation, {\it Comm. in Analysis and
Geometry, Vol. 3}, 1995, pp 523-542.

\bibitem{DD2} Daskalopoulos,P., del Pino M.A., Type II collapsing of maximal Solutions to the Ricci flow in $\R^2$, to appear in  Ann. Inst. H. Poincar Anal. Non Linaire. 

\bibitem{DH} Daskalopoulos, P., Hamilton, R., Geometric Estimates for
the Logarithmic Fast Diffusion  Equation, {\it Comm. in Analysis
and Geometry}, 2004, to appear.

\bibitem{DS} Daskalopoulos, P., Sesum, N., Eternal solutions to the Ricci flow on $\R^2$,
preprint. 

\bibitem{ERV}
Esteban, J.R., Rodr\'{\i}guez, A., Vazquez, J.L., A nonlinear heat
equation with singular diffusivity, {\em Arch. Rational Mech.
Analysis, 103}, 1988, pp. 985-1039.

\bibitem{GPV}
Galaktionov,V.A., Peletier, L.A., Vazquez, J.L.,
Asymptotics of the fast-diffusion equation with critical exponent;
{\em Siam. J. Math. Anal., 31}(1999), 1157--1174.


\bibitem{GVb}
Galaktionov, Victor A.; Vazquez, Juan Luis A stability technique
for evolution partial differential equations. A dynamical systems
approach. Progress in Nonlinear Differential Equations and their
Applications, 56. Birkhauser Boston, Inc., Boston, MA, 2004.




\bibitem{HY} Hamilton, R., Yau, S-T, The Harnack estimate for the Ricci flow on a  surface - Revisited, {Asian J. Math}, Vol 1, No 3, pp. 418-421. 


\bibitem{H1} Hamilton, R., The Ricci flow on surfaces,
{\em Contemp. Math., 71},  Amer. Math. Soc., Providence, RI, 1988,
pp 237-262.



\bibitem{H3}  Hamilton, R.  The formation of singularities in the
Ricci flow, {\em  Surveys in differential geometry, Vol. II} pp
7--136, Internat. Press, Cambridge, MA, 1995.

\bibitem{H4} Hamilton, R., The Harnack estimate for the Ricci Flow,
J. Differential Geometry {\bf 37} (1993) pp 225-243.

\bibitem{HP}
Herrero, M. and Pierre, M., The Cauchy problem for $u_t = \Delta
u^m$ when $0<m<1$, {\em Trans. Amer. Math. Soc., 291}, 1985, pp.
145-158.


\bibitem{Hsu1}
Hsu, S.-Y; Dynamics of solutions of a singular diffusion equation;
{\it  Adv. Differential Equations} {\bf 7 } (2002), no. 1, 77--97.

\bibitem{Hsu2} Hsu, Shu-Yu; Asymptotic profile of solutions of a
singular diffusion equation as $t\to\infty$, {\it Nonlinear Anal.}
{\bf 48} (2002), no. 6, Ser. A: Theory Methods, 781--790.

\bibitem{Hsu3}
Hsu, S.-Y; Large time behaviour of solutions of the Ricci flow
equation on $R\sp 2$, {\it Pacific J. Math.} {\bf 197} (2001), no.
1, 25--41.

\bibitem{Hsu4}
Hsu, S.-Y; Asymptotic behavior  of solutions of the equation
$u_t=\Delta \log u$ near the extinction time, {\em Advances in
Differential Equations}, {\bf 8}, No 2, (2003), pp 161--187.

%

\bibitem{Hui1}
Hui, K.-M. Singular limit of solutions of the equation $u\sb
t=\Delta({u\sp m}/m)$ as $m\to0$. Pacific J. Math. 187 (1999), no.
2, 297--316.


\bibitem{K} King, J.R., Self-similar behavior for the equation of fast
nonlinear diffusion, {\em phil. Trans. R. Soc., London, A 343},
(1993), pp 337--375.

\bibitem{RVE}  Rodriguez, A;  Vazquez, J. L.; Esteban, J.R, The
maximal solution of the logarithmic fast diffusion equation in two space
dimensions, {\em  Adv. Differential Equations 2} (1997), no. 6, pp
867--894.

\bibitem{W1}  Wu, L.-F.,
A new result for the porous medium equation  derived from the
Ricci flow, {\em  Bull. Amer. Math. Soc., 28}, 1993, pp 90-94.

\bibitem{W2} Wu, L.-F.,
The Ricci Flow on Complete ${\bf R}^2$, {\em Communications in
Analysis and Geometry, 1}, 1993, pp 439-472.

\end{thebibliography}
\end{document}